\documentclass[conference]{IEEEtran}
\IEEEoverridecommandlockouts
\usepackage{cite}
\usepackage{amsmath,amssymb,amsfonts}
\usepackage{mhchem}
\usepackage{gensymb}
\usepackage{algorithmic}
\usepackage[caption=false,font=normalsize,labelfont=sf,textfont=sf]{subfig}
\usepackage{graphicx}
\usepackage{textcomp}
\usepackage{xcolor}
\def\BibTeX{{\rm B\kern-.05em{\sc i\kern-.025em b}\kern-.08em
    T\kern-.1667em\lower.7ex\hbox{E}\kern-.125emX}}
\usepackage{booktabs}
\usepackage{multirow}
\usepackage{hyperref}
\usepackage{comment}
\usepackage{url}

\usepackage{breakurl}

\newcounter{subeq}

\newcommand{\figspace}{\vspace{-0.45cm}}

\makeatletter
\newcounter{savesubsection}
\newcounter{apdxsubsection}
\newcommand\unappendix{\par
  \setcounter{apdxsubsection}{\value{subsection}}%
  \setcounter{subsection}{\value{savesubsection}}%
  \gdef\thesubsection{\@arabic\c@subsection}}
\makeatother

\makeatletter
\newcommand{\settitle}{\@maketitle}
\makeatother

\begin{document}
\title{Optimization of Hybrid Power Plants: \newline When Is a Detailed Electrolyzer Model Necessary?}


\author{%
\IEEEauthorblockN{Manuel Tobias Baumhof, Enrica Raheli, Andrea Gloppen Johnsen, and Jalal Kazempour}
\IEEEauthorblockA{Department of Wind and Energy Systems, Technical University of Denmark, Kgs. Lyngby, Denmark\\ 
$\{$mtba, enrah, anglopj, jalal$\}$@dtu.dk 
}
}
\maketitle

\thispagestyle{plain}
\pagestyle{plain}

\begin{abstract}
Hybrid power plants comprising renewable power sources and electrolyzers
are envisioned to play a key role in accelerating the transition towards decarbonization. It is common in the current literature to use simplified operational models for electrolyzers. It is still an open question whether this is a good practice, and if not, when a more detailed operational model is necessary. This paper answers it by assessing the impact of adding different levels of electrolyzer details, i.e., physics and operational constraints, to the optimal dispatch problem of a hybrid power plant in the day-ahead time stage. Our focus lies on the number of operating states (on, off, standby) as well as the number of linearization segments used for approximating the non-linear hydrogen production curve. For that, we develop several mixed-integer linear models, each representing a different level of operational details.
We conduct a thorough comparative ex-post performance analysis under different price conditions, wind farm capacities, and minimum hydrogen demand requirements, and discuss under which operational circumstances a detailed model is necessary. In particular, we provide a case under which a simplified model, compared to a detailed one, results in a decrease in profit of 1.8\%  and hydrogen production of 13.5\% over a year.
The key lesson learned is that a detailed model potentially earns a higher profit  in circumstances under which the electrolyzer operates with partial loading. This could be the case for  a certain range of electricity and hydrogen prices, or limited wind power availability. The detailed model also provides a better estimation of true hydrogen production, facilitating the logistics required.
\end{abstract}

\begin{IEEEkeywords}
hybrid power plants, electrolyzer, hydrogen, mixed-integer linear programming
\end{IEEEkeywords}
\linespread{0.89}
\section{Introduction}
\subsection{Background}
In order to limit global warming to a maximum of 1.5~$\degree$C, greenhouse gas emissions must be reduced to net zero by 2050, as called for in the European Green Deal 2019 \cite{EUGreenDeal}. Renewable hydrogen produced through electrolysis could aid in two major challenges on the path towards the net zero goal. First, electrolyzers can act as flexible loads and therefore potential frequency restoration ancillary service providers, contributing to maintaining the power balance in power systems with increased penetration of renewable energy sources. Second, renewable hydrogen can be further synthesized into other green fuels, 
eventually enabling decarbonization in the hard-to-abate sectors, such as heavy transport and industry. 

Hybrid power plants comprising of renewable power sources (wind and/or solar) and electrolyzers are the key components to accelerate the current energy transition through hydrogen \cite{PtX_strategy_EU}. 
Nonetheless, uncertainties in terms of the cost-benefit of electrolyzers in the long run have challenged the widespread investment in said technologies and thereby large-scale production of renewable-based green hydrogen \cite{iea_hydrogen_review}. In Denmark, there is currently a special focus on green hydrogen at the governmental level and also, among the regulator, system operator, and many industry stakeholders, envisioning a large deployment of electrolyzers and other power-to-X facilities in the coming years. In 2021 the Danish government published a strategy for the national power-to-X development, aiming to build 4 to 6 GW of electrolysis capacity by 2030, doubling the current Danish peak demand \cite{PtX_strategy_DK}. This emerging trend is not limited to Denmark, and many other countries both in Europe and globally see hydrogen as a key solution for the realization of green societies of the future \cite{PtX_strategy_EU,India}. 



\subsection{Aim and Literature Review}
It is a common practice in the current literature to use a simplified operational model for electrolyzers e.g., by using a constant power-to-hydrogen conversion ratio irrespective of whether the electrolyzer operates in full capacity or not \cite{Mancarella,MATUTE201917431,MATUTE20211449,martin}.
In addition, some papers do not consider operational states of the electrolyzer \cite{Mancarella,martin}.
This paper challenges these simplification practices. While a simplified model works satisfactorily under certain operational circumstances, there are several other circumstances under which a simplified one yields a sub-optimal operation of electrolyzers, underestimating their value. This paper answers  when a detailed operational model should be applied, and to what extent the profit and hydrogen production can be increased by using a detailed model. We will also discuss to what extent a detailed model brings additional computational burden. 

In general, two main physical aspects of electrolyzers need to be modeled for operation in the day-ahead time stage:
\begin{enumerate}
    \item \textit{Electrolyzer efficiency}: The power-to-hydrogen conversion efficiency is a function of the power consumption of the electrolyzer. To accurately model the hydrogen production of the electrolyzer, the varying efficiency should be captured, which introduces non-linearities to the model. The simple models usually use a constant efficiency, while more accurate modeling incorporates the non-linearities, which can be later linearized.
    
    \item \textit{Number of operating states}: Proper operational modeling of electrolyzers may require introducing three states, namely on, off, and standby, to ensure no hydrogen production below a given minimum allowed partial loading, for which additional binary variables are needed. Many papers in the literature do not even model states, thus assuming the electrolyzer is always on, or model two states only, i.e., on and off, similar to conventional power generators\footnote{We will discuss later in Section \ref{results} that under some operational conditions, a two-state model including on and standby states works well too.}. 
\end{enumerate}

Various studies have incorporated different levels of operational details of the electrolyzer into their optimization problems. In \cite{MATUTE201917431} and \cite{MATUTE20211449}, a constant efficiency is applied but two and three states are modeled, respectively, by adding binary variables. In \cite{VARELA20219303}, three states are modeled, while assuming a linear hydrogen production curve, despite showing that the production curve is not well approximated by a first-order interpolation. 
A hybrid power plant including an electrolyzer is modeled in \cite{pavic2022}, where the non-linear hydrogen production is linearized between two points, with a single binary variable representing the on/off state of the electrolyzer. In \cite{beerbuhl2014} a quadratic production curve is applied and the resulting non-linear program is eventually solved by a heuristic.
In \cite{ShiYou}, three states are included, and differently from the other papers, the operating temperature is considered as a variable, providing an extra degree of freedom in the electrolyzer operation. This model allows to take into account the temperature impact on the conversion efficiency and the quality of the generated heat. The non-linear hydrogen production is then linearized around a fixed reference operating point to formulate the problem as a mixed-integer linear programming (MILP) problem.

\subsection{Contributions and Paper Organization}
To the best of our knowledge, there is a lack of a comprehensive analysis in the current literature, identifying the operational circumstances under which a simple model ends up in a sub-optimal operation of electrolyzers, resulting in a reduced profit and hydrogen production\footnote{Reference \cite{ShiYou} provides a similar analysis, however, the Faraday efficiency is assumed to be one. The consequences of this assumption will be further discussed in Section \ref{eff_prod}.}.
This paper bridges such a gap through the following contributions:
\begin{itemize}
\item To embed constraints describing the physics of electrolyzers while keeping the final model as a MILP, 
\item To thoroughly investigate ex-post the impact of the inclusion of different operational details on the final profit of the hybrid power plant and  the amount of hydrogen produced, 
    \item and finally, to provide a set of recommendations  in terms of including operational details of electrolyzers, depending on the application, the range of electricity prices, and the hydrogen price.
\end{itemize}

Without loss of generality, this  paper focuses on alkaline electrolyzers,
as they are currently the most mature technology \cite{gotz2016}. The proposed model can be extended to other low-temperature electrolyzers, such as polymer electrolyte membrane (PEM). More operational characteristics may be necessary for modeling 
solid-oxide electrolyzers (SOEC).

The rest of the paper is organized as follows. Section \ref{physics} describes  the electrolyzer physics, focusing on the operating states and the hydrogen production curve. Section \ref{model} provides the proposed MILP, representing all three states of the electrolyzer. Section \ref{results} discusses the impact of the electrolyzer modeling choices by means of a test case and a thorough sensitivity analysis. Section \ref{conclusion} concludes the paper. In the Appendix, the day-ahead price range where electrolyzer details matter is analytically formulated. Finally, the Online Companion \cite{Online_Companion} provide two MILPs (simpler than the one proposed in Section \ref{model}), both representing two states of the electrolyzer only, where one is a model with on-off states, and the other one is a model with on-standby states.

\section{Electrolyzer physics}
\label{physics}
The core of the renewable-hydrogen hybrid power plant is the electrolyzer, where water is decomposed into hydrogen and oxygen by means of electrical power. The physics and operating characteristics of alkaline electrolyzers are described in this section and will be 
formulated as a set of mixed-integer linear constraints in Section \ref{model}. 

\subsection{States}
To describe and model the real operation of an alkaline electrolyzer, it is necessary to distinguish three different states:
\subsubsection{On state} the electrolyzer operates within its feasible load range, consuming power and producing hydrogen with a conversion efficiency that depends on the partial load, which will be explained in Section \ref{eff_prod}. 
The minimum operating power for alkaline electrolyzers is around 15-20\% of the nominal power, below which the electrolyzer must go into standby or off. 
\subsubsection{Standby state} the electrolyzer does not produce any hydrogen but consumes the power needed to maintain the system temperature and pressure so that it can rapidly resume production. The value of the standby power consumption is not usually disclosed by manufacturers, but values between 1-5\% of the electrolyzer full load capacity have been adopted in the literature \cite{VARELA20219303, MATUTE20211449, MATUTE201917431}. The time needed to switch from standby to on, i.e., a warm start-up is of the order of 30 seconds \cite{ MATUTE20211449}.
\subsubsection{Off state} the electrolyzer is shut down completely and does not consume any power nor produce any hydrogen. However, to switch back to on, a significant amount of electricity is needed, corresponding to a cold start-up cost. Moreover, at least 20 minutes are necessary before resuming hydrogen production \cite{ MATUTE20211449}. Apart from the introduced cold start-up cost and start-up time, the frequent shut down of the electrolyzer may have a negative impact on the device degradation and lifetime \cite{URSUA201612852}.

\subsection{Efficiency and Production Curve}
\label{eff_prod}
The conversion efficiency of electricity into hydrogen is not constant but depends on the \textit{partial load}, i.e., the ratio between power consumption at a specific time and the nominal power of the electrolyzer. The variation of the efficiency based on the operating set-point is mainly due to two phenomena: (\textit{i}) the current-voltage relationship, also called the polarization curve, and (\textit{ii}) the Faraday efficiency. We explain both phenomena in the following.

The current-voltage relationship describes the voltage increase (also called over-voltage or over-potential) with increasing current density, due to different losses, as explained in \cite{SANCHEZ20203916} and \cite{ShiYou}. Ulleberg \cite{Ulleberg} introduced a widely adopted empirical formulation that describes the relationships between voltage, current density, and electrolyzer operating temperature. To further take into account the operating pressure, this formulation was modified by Sanchez et al. \cite{Sanchez}. For a given temperature and pressure, this can be formulated as
\begin{align}
    U^{\rm{cell}}(i)  = U^{\rm{rev}} + K_1 i + K_2 \text{log} (K_3 i +1),
    \label{eq:voltage}
\end{align}
where $U^{\rm{cell}}(i)$ is the cell voltage as a function of  the current density $i$. In addition, $U^{\rm{rev}}$ is the open-circuit voltage (i.e., voltage corresponding to current density equal to zero). The parameters $K_1$, $K_2$, $K_3$ are constants obtained from experimental data and can be found in \cite{Sanchez}. Voltage $U^{\rm{rev}}$ can be calculated for a specific operating temperature according to an empirical equation that can be found in \cite{Sanchez}. The power consumed by the electrolyzer $p^{\rm{e}}(i)$ can be calculated as
\begin{align}
    & p^{\rm{e}}(i) = U^{\rm{cell}}(i) i A, \label{eq:power}
\end{align}
where $A$ is the total area of the cells composing the electrolyzer. 
\begin{figure}[t]
\centering
\includegraphics[width=\linewidth]{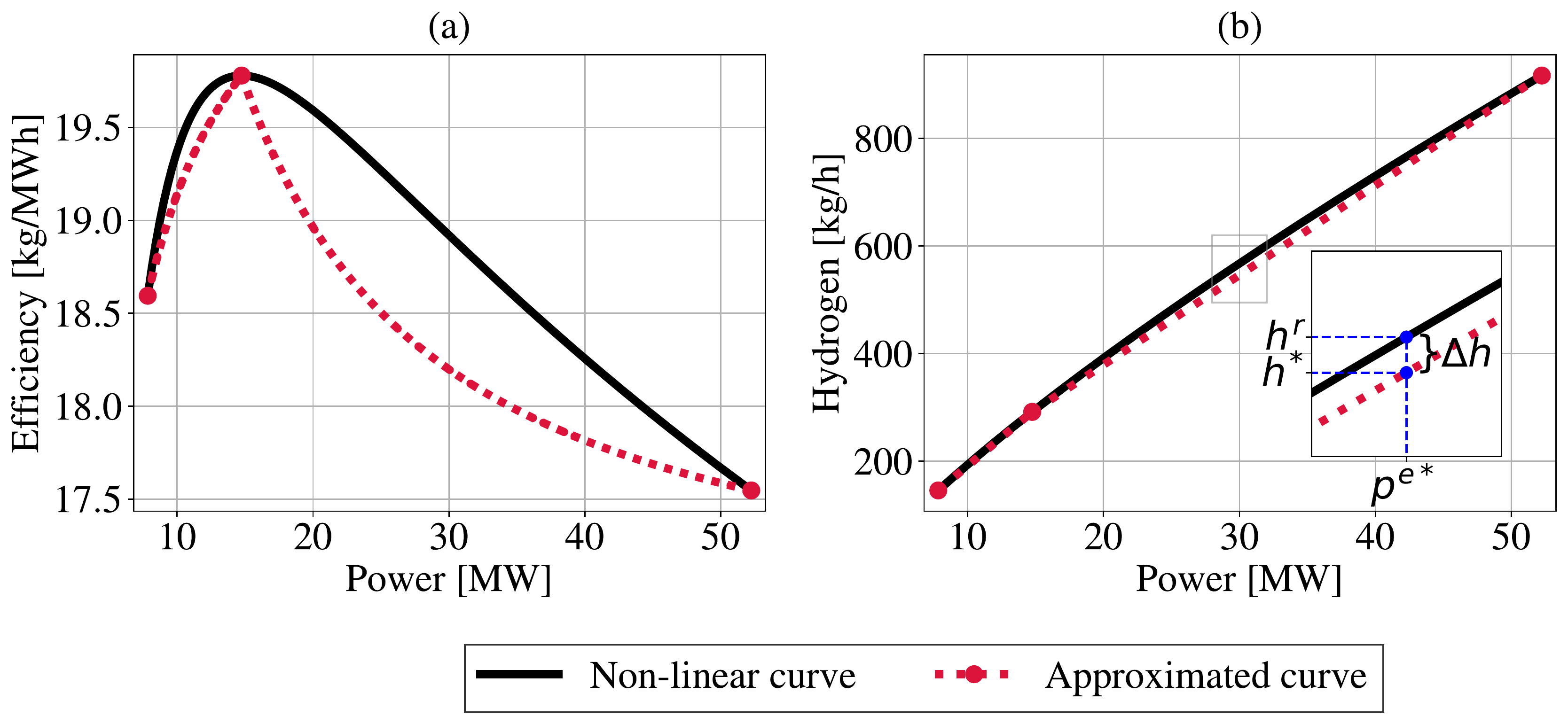}%
\caption{Plot (a): the efficiency curve, and plot (b): the hydrogen production curve of a 52.25-MW alkaline electrolyzer, as a function of the electric power consumption, working at 90 \degree C and 30 bar. The black curves represent the original non-linear curves. Approximated by two segments, the red curve in plot (b) is the piecewise linearized hydrogen production curve. The non-linear efficiency curve corresponding to this piecewise linearization is represented by the red curve in plot (a).  In our formulation, we will only use the red piecewise linear production curve in plot (b). The inner plot of (b) shows the hydrogen production discrepancy $\Delta h$ between original and approximated curves, for a given power consumption level. }
\label{fig:fig_eff_prod}
\figspace
\end{figure}
The Faraday law calculates the hydrogen production $h(i)$ of the electrolyzer as
\begin{align}
    h(i) = 3600\cdot \frac{\eta^{\rm{F}}(i) M^{\rm{H_2}} i A}{2F}, \label{eq:h_prod}
\end{align}
where $h(i)$ is the hydrogen production rate in kg/h, $M^{\rm{H_2}}$ is the molar mass of hydrogen in kg/mol, $F$ is the Faraday constant, and $\eta^{\rm{F}}(i)$ is the Faraday efficiency as a function of current density. 
The latter is defined as the ratio between the actual and the theoretical maximum amount of hydrogen produced. The difference between actual and theoretical output is explained in \cite{Ulleberg}, and it increases significantly when the electrolyzer is working at low-current densities. 
In \cite{Sanchez}, an empirical expression that captures the relationship between the Faraday efficiency and the current density at a given temperature is provided: $\eta^{\rm{F}}(i)$ is close to one for higher current densities, and it drops to zero when reducing the current. 
The electrolyzer efficiency is defined as
\begin{align}
    \eta(i) = \frac{h(i)}{{p^{\rm{e}}}(i)}, \label{eq:eff}
\end{align}
where generally $\eta(i)$ is expressed in kg/MWh. For different values of $i$, the black curve in Fig.~\ref{fig:fig_eff_prod}(a) shows efficiency $\eta(i)$ versus power consumption ${p^{\rm{e}}}(i)$. In addition, the black curve in Fig.~\ref{fig:fig_eff_prod}(b) shows the hydorgen production $h(i)$ versus power consumption ${p^{\rm{e}}}(i)$. For notational clarity, we drop $(i)$ in the rest of the paper. The black curves in  Fig.~\ref{fig:fig_eff_prod} show that the model is non-linear. The efficiency has a peak at around 30\% of the load. This characteristic peak in the efficiency curve is not captured when a constant conversion efficiency is used, as done in \cite{VARELA20219303, MATUTE20211449, Mancarella}, 
or when the Faraday efficiency is assumed to be equal to one in the entire feasible operating range, as done in \cite{ShiYou}. 

To keep the final problem a MILP, but describe the hydrogen production with more details, we use a piecewise linearization of the hydrogen production curve as shown by the red curve in  Fig.~\ref{fig:fig_eff_prod}(b), for two linearization segments. For each segment $s \in \mathcal{S}$, the $A_s$ (slope) and $B_s$ (intercept) coefficients of the line can be calculated such that the approximated hydrogen production is $A_s p^{\rm{e}} + B_s$. Later we will define a binary variable indicating which segment is active. The proposed approximation is exact only at the segment endpoints (i.e., linearization points), otherwise, it is an underestimation of the original non-linear curve. For example, the optimal  power set-point $p^{\rm{e^*}}$ in the inset of Fig.~\ref{fig:fig_eff_prod}(b) corresponds to the hydrogen production $h^{*}$ according to the proposed piecewise linear model with two segments\footnote{Symbol ${}^*$ refers to the optimal value.}. However, the actual hydrogen realization based on the electrolyzer physics is $h^{\rm{r}}$. The hydrogen production difference  $\Delta h$ is reduced by increasing the number of segments, and the effect of the hydrogen surplus obtained when choosing only one segment, as done in \cite{VARELA20219303}, is discussed in Section \ref{results}.

According to this piecewise linear formulation for the hydrogen production curve, the efficiency $\eta$ for segment $s$ can be calculated based on \eqref{eq:eff}, resulting in $\eta = A_s + \frac{B_s}{{p^{\rm{e}}}}$. This is depicted by the red dotted curve in Fig.~\ref{fig:fig_eff_prod}(a), given two linearization segments used. Note that it does not present a linear behavior. However, this non-linear efficiency curve does not appear in our optimization problem. The hydrogen production curve is used instead, which is linearized through segments, as illustrated by the red dotted curve in Fig.~\ref{fig:fig_eff_prod}(b).


\section{Problem formulation}
\label{model}

\begin{figure}[t]
\centering
\includegraphics[width=\linewidth]{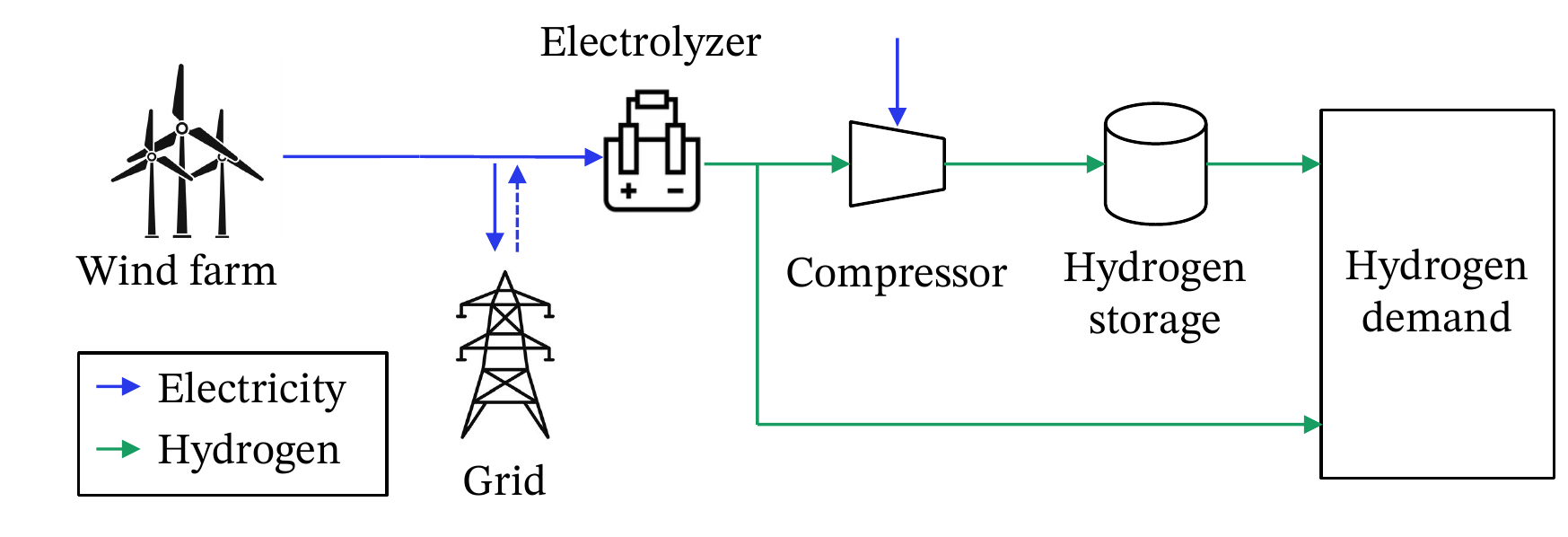}%
\caption{Schematic representation of a hybrid power plant.}
\label{fig:hpp_schematic}
\figspace
\end{figure}

We consider a hybrid power plant, as depicted in  Fig.~\ref{fig:hpp_schematic}, consisting of a wind farm, an electrolyzer, a hydrogen compressor, and a hydrogen storage. The generated wind power can be either sold to the grid at the electricity market price, or consumed by the electrolyzer to produce 100\% renewable-based green hydrogen.
The hydrogen produced can either be directly delivered to the demand or temporarily stored in an on-site hydrogen storage, with an associated cost for compressing the gas. The dashed blue line in  Fig.~\ref{fig:hpp_schematic} represents the option to buy electricity from the grid only to supply the electrolyzer's standby power when there is no wind power.

The hydrogen price is assumed to be a single-value constant, and the hybrid power plant serves a minimum daily hydrogen demand. 
We assume the plant has perfect foresight of future wind power production and electricity price. Given the 1-hour time resolution in our model, we neglect the ramping limitation which are typically around $\pm$20\% of the nominal power per second \cite{VARELA20219303}, as well as the warm and cold start-up times of the electrolyzer.

For the optimal operation of the hybrid power plant, we develop a complete MILP in Section \ref{ssec:model_3states} accounting for three states of the electrolyzer and then provide two simplified counterparts in Section \ref{ssec:model_2states}, each with two states of the electrolyzer. 

\textit{Notation}: All parameters are upper-case or Greek letters, whereas all variables are lower-case letters.  All binary variables are noted by \textit{z}.




\subsection{Three-state Model}
\label{ssec:model_3states}
The most complete MILP includes the objective function \eqref{eq:2_objective} constrained by \eqref{eq:market2}-\eqref{last_con3x}. 
\subsubsection{Objective function}
%
%
Over the set of hours $\textit{t} \in \mathcal{T}$, the objective function \eqref{eq:2_objective} maximizes the total profit of the hybrid power plant as
\begin{align}
    & \underset{\bf{x}}{\rm{max}} \sum_{t \in \mathcal{T}} p_t \lambda^{\rm{DA}}_t + d_t \lambda^{\rm{h}} - p^{\rm{in}}_t \lambda^{\rm{in}}_t - z^{\rm{su}}_t \lambda^{\rm{su}},
   \label{eq:2_objective}
\end{align}
where the variable set $\bf{x}$ will be defined later. The first term corresponds to selling power $p_t$ to the grid at the day-ahead electricity market price $\lambda^{\rm{DA}}_t$. The second term pertains to delivered hydrogen $d_t$ at a fixed price $\lambda^{\rm{h}}$. 
%
The third term represents the cost for purchasing standby power $p^{\rm{in}}_t$ to support the electrolyzer's standby state in case the wind power is insufficient. The corresponding price is $\lambda^{\rm{in}}_t = \lambda^{\rm{DA}}_t + \lambda^{\rm{TSO}}$, where $\lambda^{\rm{TSO}}$ is the grid tariff imposed by the Transmission System Operator (TSO). Finally, the fourth term corresponds to the cold start-up cost of the electrolyzer, where the binary variable $z^{\rm{su}}_t$ indicates the start-up at hour $t$,  associated with the cost per startup $\lambda^{\rm{su}}$. 

\subsubsection{Power balance}

In every hour $t$, the power $p_t$ sold in the day-ahead market is equal to the wind farm power production $P^{\rm{w}}_t$ plus power $p^{\rm{in}}_t$ bought from the grid to support the standby state of the electrolyzer, subtracted by the power consumption $p^{\rm{e}}_t$ of the electrolyzer and the power consumption $p^{\rm{c}}_t$ of the compressor, such that 
\begin{align}
     & p_t = P^{\rm{w}}_t + p^{\rm{in}}_t- p^{\rm{e}}_t - p^{\rm{c}}_t & \forall~ & t \in \mathcal{T}. \label{eq:market2} 
\end{align}

\subsubsection{Limit on $p^{\rm{in}}_t$}
The input power $p^{\rm{in}}_t$ is limited by the standby state consumption of the electrolyzer, implying that power cannot be bought from the grid to produce hydrogen: 
%
\begin{align}
    & p^{\rm{in}}_t \leq P^{\rm{sb}} z^{\rm{sb}}_t & \forall~ & t \in \mathcal{T}, \label{eq:market_in} 
\end{align}
where the parameter $P^{\rm{sb}}$ is the standby consumption, and the binary variable $z^{\rm{sb}}_t$ indicates whether the electrolyzer is in the standby mode in hour $t$.

\subsubsection{Electrolyzer operational states}
Constraint \eqref{eq:status} ensures that the electrolyzer can take only one out of three states at any hour \textit{t}, namely online, standby, or off: 
\begin{align}
    & z^{\rm{on}}_t + z^{\rm{off}}_t + z^{\rm{sb}}_t = 1 & \forall~ & t \in \mathcal{T},\label{eq:status}
\end{align}
where similar to $z^{\rm{sb}}_t$, binary variables $z^{\rm{on}}_t$ and $z^{\rm{off}}_t$ indicate whether in hour $t$ the electrolyzer is on and off, respectively.
The states are activated based on the electricity consumption of the electrolyzer. In the online state, the electricity consumption $p^{\rm{e}}_t$ of the electrolyzer can neither exceed the capacity  $C^{\rm{e}}$ nor go below a minimum load limit $P^{\rm{min}}$. In the standby state, the electricity consumption must be equal to the standby power consumption $P^{sb}$. These constraints are enforced by
\begin{align}
    & p^{\rm{e}}_t  \leq C^{\rm{e}} z^{\rm{on}}_t + P^{\rm{sb}} z^{\rm{sb}}_t & \forall~ & t \in \mathcal{T} \label{eq:elec_max}, \\
    & p^{\rm{e}}_t \geq P^{\rm{min}} z^{\rm{on}}_t + P^{\rm{sb}} z^{\rm{sb}}_t & \forall~ & t \in \mathcal{T}. \label{eq:elec_min}
\end{align}
%
To represent the cold start-up of the electrolyzer, the binary variable $z^{\rm{su}}_t$ is defined, taking the value 1 in the case of a transition from off to on state in hour \textit{t}, as enforces by constraints \eqref{eq:start} and  \eqref{eq:start_1}. Further, constraint \eqref{eq:off_to_sb}   
 ensures that the transition from an off-state to a standby-state is not allowed, to avoid bypassing of the start-up cost. 
\begin{align}
    & z^{\rm{su}}_t \geq z^{\rm{on}}_t - z^{\rm{on}}_{t-1} - z^{\rm{sb}}_{t-1} & \forall~ & t \in \mathcal{T} \setminus \{1\}, \label{eq:start} \\
    & z^{\rm{su}}_{t=1} = 0, \label{eq:start_1} \\
    & z^{\rm{off}}_{t-1} + z^{\rm{sb}}_t \leq 1 & \forall~ & t \in \mathcal{T}  \setminus \{1\}. \label{eq:off_to_sb}
\end{align}

\subsubsection{Electrolyzer hydrogen production}
The hydrogen production $h_t$ is a function of the electricity consumption of the electrolyzer. As explained in Section \ref{eff_prod}, for each segment $s \in \mathcal{S}$, a linear function of the segment power consumption $\hat{p}^{\rm{e}}_{ts}$ with slope $A_s$ and intercept $B_s$ is defined, such that
\begin{align}
        & h_t = \sum_{s \in \mathcal{S}} (A_s \hat{p}^{\rm{e}}_{ts} + B_s z^{\rm{h}}_{ts})& \forall~ & t \in \mathcal{T}, \label{eq:hy}
\end{align}
where the binary variable $z^{\rm{h}}_{ts}$ defines which segment $s$ is active in hour $t$.
Each segment is valid within a pre-defined interval of upper $\overline{P}_s$ and lower $\underline{P}_s$ power consumption levels, i.e.,
\begin{align}
    & \underline{P}_s  z^{\rm{h}}_{ts} \leq \hat{p}^{\rm{e}}_{ts} \leq  \overline{P}_s  z^{\rm{h}}_{ts}& \forall~  t \in \mathcal{T}, s \in \mathcal{S}. \label{eq:s_min}
\end{align}

Constraint \eqref{eq:prod} ensures that hydrogen production happens in the online state only, while one segment only can be active at any hour \textit{t}. In addition, \eqref{eq:power_tot} computes the total power consumption of the electrolyzer:
\begin{align}
    & z^{\rm{on}}_t = \sum_{s \in \mathcal{S}} z^{\rm{h}}_{ts} & \forall~ & t \in \mathcal{T}, \label{eq:prod}\\
    & p^{\rm{e}}_t = \sum_{s \in \mathcal{S}} \hat{p}^{\rm{e}}_{ts} + P^{\rm{sb}} z^{\rm{sb}}_t& \forall~ & t \in \mathcal{T}. \label{eq:power_tot}
\end{align}

\subsubsection{Hydrogen storage}
Constraints \eqref{eq:stor_1}-\eqref{eq:stor_2} represent the storage operation:
\begin{align}
    & h_t = h^{\rm{d}}_t + s^{\rm{in}}_t & \forall~ & t \in \mathcal{T}, \ \label{eq:stor_1}\\
    & d_t = h^{\rm{d}}_t + s^{\rm{out}}_t & \forall~ & t \in \mathcal{T}, \label{eq:stor_3}\\
     & s^{\rm{out}}_t \leq S^{\rm{out}} & \forall~ & t \in \mathcal{T}, \label{eq:stor_4}\\
      & p^{\rm{c}}_t = K^{\rm{c}} s^{\rm{in}}_t & \forall~ & t \in \mathcal{T},\label{eq:stor_5}\\
          & s_{t=1} = S^{\rm{ini}} + s^{\rm{in}}_{t=1} - s^{\rm{out}}_{t=1} &  \label{eq:stor_6}\\
    & s_t = s_{t-1} + s^{\rm{in}}_t - s^{\rm{out}}_t & \forall~ & t \in \mathcal{T} \setminus \{1\},\label{eq:stor_7}  \\
    & s_t \leq C^{\rm{s}} & \forall~ & t \in \mathcal{T}. \label{eq:stor_2}
\end{align}
The hydrogen produced $h_t$ can either go directly to the demand $h^{\rm{d}}_t$ or be injected into the hydrogen storage $s^{\rm{in}}_t$, as enforced by \eqref{eq:stor_1}. The total hydrogen $d_t$ delivered to the demand is equal to the sum of hydrogen directly from the electrolyzer and that from the storage $s^{\rm{out}}_t$, as per \eqref{eq:stor_3}. The storage output of every hour is limited by the output flow capacity $S^{\rm{out}}$ in \eqref{eq:stor_4}. Further, the compressor consumes power $p^{\rm{c}}$ to compress the hydrogen injected into the storage. Assuming adiabatic compression, the compression coefficient $K^{\rm{c}}$ can be calculated, as proposed by \cite{ShiYou}. The power consumption for compression is then \eqref{eq:stor_5}. The state of charge of the hydrogen storage in the initial and following hours is calculated by \eqref{eq:stor_6} and \eqref{eq:stor_7}, where $S^{\rm{ini}}$ is the hydrogen initially stored in the storage at the beginning of time horizon $\mathcal{T}$. The storage hydrogen mass capacity $C^{\rm{s}}$ is enforced by \eqref{eq:stor_2}.
Note that we do not impose any constraint for the energy stored at the end of time horizon $\mathcal{T}$. Therefore, pursuing profit maximization in this time horizon, the hybrid power plant will leave the storage empty in the last hour\footnote{One can enforce a constraint on the minimum stored hydrogen at the end of the time horizon, or add a value for this stored energy to the objective function.}.

\subsubsection{Hydrogen demand} Imagine within the underlying time horizon $\mathcal{T}$, which could be, for example, a year, there are $N$ number of time subsets, e.g., 365 days, indexed by $n$, such that there is a minimum hydrogen demand for each $n$:
\begin{align}
    & \sum_{t \in \mathcal{H}_n} d_{t} \geq D_n^{\rm{min}} & \forall~ & n \in \{1, ..., N\}, \label{eq:demand}
\end{align}
where $\mathcal{H}_n$ is the set of hours within time subset $n$.


\subsubsection{Variable declaration} Constraint \eqref{last_con1} declares the non-negativity conditions:
\begin{align}
        &   d_t,  h_t, h^{\rm{d}}_t, p_t, p^{\rm{c}}_t, p^{\rm{in}}_t, \hat{p}^{\rm{e}}_{ts},   s_t, s^{\rm{in}}_t, s^{\rm{out}}_t \in \mathbb{R}^{+}. \ \label{last_con1}
\end{align}

Constraint \eqref{last_con2} lists binary variables:
\begin{align}
        & z^{\rm{su}}_t, z^{\rm{h}}_{ts}, z^{\rm{on}}_t, z^{\rm{off}}_t, z^{\rm{sb}}_t \in \{0,1\}. \label{last_con2}
\end{align}

Therefore, the total number of binary variables is $|\mathcal{T}|(4+|\mathcal{S}|)$ binaries, where $|\mathcal{T}|$ and $|\mathcal{S}|$, respectively, are the number of hours and the number of segments used to linearize the hydrogen production curve.
Finally, the variable set $\bf{x}$ is defined as
\begin{align}
         \bf{x} = \{&d_t,  h_t, h^{\rm{d}}_t, p_t, p^{\rm{c}}_t, p^{\rm{in}}_t, \hat{p}^{\rm{e}}_{ts}, \notag \\ 
        &s^{\rm{in}}_t, s_t,  s^{\rm{out}}_t, z^{\rm{su}}_t, z^{\rm{h}}_{ts}, z^{\rm{on}}_t, z^{\rm{off}}_t, z^{\rm{sb}}_t \}. \label{last_con3x}
\end{align}
Accordingly, in addition to $|\mathcal{T}|(4+|\mathcal{S}|)$ number of binary variables, we have $|\mathcal{T}|(9+|\mathcal{S}|)$ number of continuous variables.

\subsection{Two-state Models}
\label{ssec:model_2states}
The optimal operation problem \eqref{eq:2_objective}-\eqref{last_con3x} of the hybrid power plant accounting for three states of the electrolyzer  can be simplified if two states only are considered, either  on-off states or on-standby states. Both result in MILPs. 

In the latter, i.e., the MILP with on-off states, one binary variable (instead of three) per hour $t$ is sufficient, such that it indicates whether the electrolyzer in the given hour is on or off. The resulting MILP is provided in \cite{Online_Companion}. 
The total number of binary variables in this MILP is $|\mathcal{T}|(2+|\mathcal{S}|)$.

Similarly, a single binary variable per hour $t$ is enough in the MILP with on-standby states, indicating whether the electrolyzer is online or in standby mode. Also,
the start-up binary variable is not needed. The corresponding MILP is given in \cite{Online_Companion}, where among three MILPs, we need the lowest number of binary variables, i.e.,  $|\mathcal{T}|(1+|\mathcal{S}|)$. 

\section{Numerical Study}
\label{results}
We apply the proposed MILPs of Section \ref{model} to a case study and investigate how the optimal operation of the hybrid power and the resulting profit change by adding more operational details of the electrolyzer. All source codes and input data are publicly shared\footnote{GitHub: \url{https://github.com/mtba-dtu/detailed-electrolyzer-model}}. We consider several options for the number of linearization segments, i.e., $|\mathcal{S}|$, used to approximate the hydrogen production curve of the electrolyzer, including 1, 2, 4, 8, and 12 segments. Also, we consider three options for the number of electrolyzer states: three states on-off-standby ($\rm{OOS}$), two states on-standby ($\rm{OS}$), and two states on-off ($\rm{OO}$). In the rest of this section, we will refer to various models as, for example, $\rm{OOS}$-12, implying we consider three states ($\rm{OOS}$) with 12 segments. 
Finally, we conduct a sensitivity analysis to explore the impact of various input parameters, such as wind farm capacity, hydrogen demand, and hydrogen price, on the operation of the hybrid power plant.

\subsection{Case Study}
\label{ssec:case}
We consider a hybrid power plant whose structure equals the one in  Fig.~\ref{fig:hpp_schematic}, and its input data is provided in Tab.~\ref{tab:case_details}. 
The capacity of the wind farm  is 104.5 MW, corresponding to 11 V164-9.5 MW™ Vestas turbines, located in Køge Bay, Denmark.  
The electrolyzer capacity is set to 50\% of the wind farm capacity, amounting to 52.25 MW. The modeling horizon spans one year with an hourly temporal resolution. We apply hourly electricity price data for 2019, as price data for the following years might be distorted by macroeconomic impacts, such as COVID-19. Day-ahead electricity prices for the East Denmark area (DK2) are obtained from ENTSO-e Transparency platform \cite{entso-e} and hourly historical wind capacity factors at the given location for 2019 are retrieved from the Renewable.ninja web platform \cite{renewables.ninja_2}. The average yearly capacity factor for the selected location is 43.7\%. The hybrid power plant is only allowed to buy power from the grid to keep the electrolyzer in standby mode, in case the wind power is insufficient. In that case, the electricity is bought at the hourly day-ahead market price plus the grid tariff of the TSO. Since the wind farm is located in DK2, the consumption tariff imposed by the Danish TSO, Energinet, is applied\footnote{Source: \url{https://energinet.dk/El/Elmarkedet/Tariffer/Aktuelle-tariffer/}}.
The minimum daily demand can be met by the full-load operation of the electrolyzer for around four hours. The hydrogen storage is scaled to store all hydrogen produced if the electrolyzer operates at full capacity for 24 consecutive hours.


\begin{table}[t]
    \centering
    \caption{Input data for the case study}
    \label{tab:case_details}
    \begin{tabular}{@{}lllrl@{}}
    \toprule
    Wind farm & Capacity & $C^{\rm{w}}$  & 104.5 & MW \\ \midrule
    \multirow{8}{*}{Electrolyzer} & Capacity & $C^{\rm{e}}$ & 52.25 & MW  \\
     & Standby load & $P^{\rm{sb}}$ & 0.52 & MW \\
     & Minimum load & $P^{\rm{min}}$ & 7.84 & MW  \\
     & Pressure &  & 30 & bar \\
     & Temperature &  & 90 & °C \\
     & Max. current density &  & 5,000 & A/m$^2$ \\
     & Start-up cost & $\lambda^{\rm{su}}$ & 2,612.50 & € \cite{VARELA20219303}\\
     & TSO tariff & $\lambda^{\rm{TSO}}$ & 15.06 & €/MWh \\ \midrule
    \multirow{2}{*}{Storage} & Capacity & $C^{\rm{s}}$ & 22,000 & kg \\
     & Maximum output & $S^{\rm{out}}$ & 912.13 & kg/h \\ \midrule
    Compressor & Consumption coefficient & $K^{\rm{c}}$ & $0.0012$ & MWh/kg \\ \midrule
    \multirow{2}{*}{Hydrogen} & Price & $\lambda^{\rm{h}}$ & 2.10 & €/kg \\
     & Minimum demand & $D_n^{\rm{min}}$ & 3,667 & kg/day \\ \bottomrule
    \end{tabular}
    \figspace
\end{table}

\subsection{Impact of the Number of Segments}  
\label{ssec:segments}
Let us consider the $\rm{OOS}$ case with three states, for which we solve the proposed MILP \eqref{eq:2_objective}-\eqref{last_con3x}. 
We start with $\rm{OOS}$-1, where $|\mathcal{S}| = 1$. This means the original non-linear hydrogen production curve, depicted in  Fig.~\ref{fig:fig_eff_prod}(b), is approximated by a single linear curve. Here, the minimum power consumption $P^{\rm{min}}$ and the capacity $C^{\rm{e}}$ of the electrolyzer are taken as two endpoints. By moving to $\rm{OOS}$-2, where  the number of segments $|\mathcal{S}|$ is 2, we consider an additional point $P^{\rm{\eta,max}}$, which refers to the power consumption level corresponding to the peak in the efficiency curve in  Fig.~\ref{fig:fig_eff_prod}(a).
%
By increasing $|\mathcal{S}|$ to 4, and then to 8, the mean load value between existing points is added, splitting one segment into two. The same procedure but only on the right side of $P^{\rm{\eta,max}}$ is applied when we move from $\rm{OOS}$-8 to $\rm{OOS}$-12, as this side covers over around 70\% of the feasible operating range. With the adoption of this procedure, all cases from $\rm{OOS}$-2 to $\rm{OOS}$-12 include the point $P^{\rm{\eta,max}}$. In addition, points are not removed when refining the discretization. By adding more segments, the hydrogen production curve and thus the electrolyzer efficiency with partial loading is more accurately represented.

\begin{figure}[t]
\centering
\includegraphics[width=\linewidth]{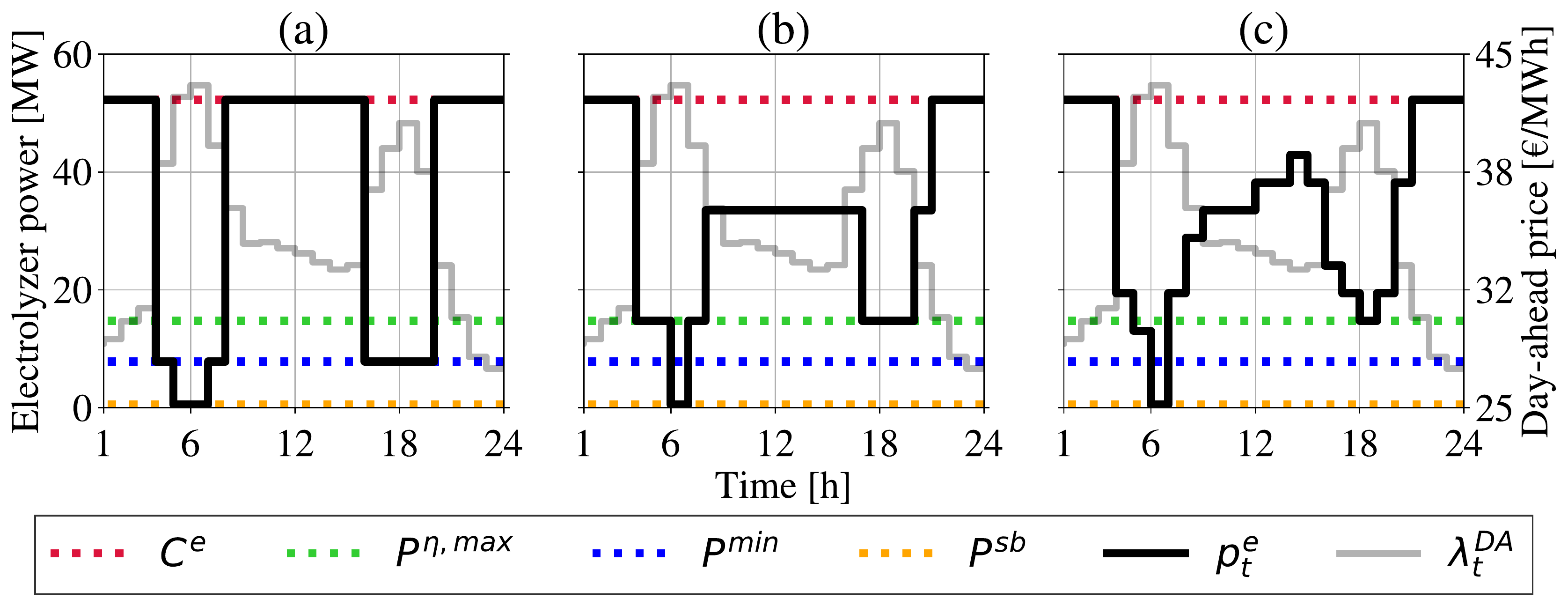}%
\caption{The power consumption schedule of the electrolyzer ($p^{\rm{e}}_t$) in an example high-wind day when its hydrogen production curve is linearized by (a) 1,  (b) 4, and (c) 12 segments. These three plots, from left to right, correspond to cases $\rm{OOS}$-1, $\rm{OOS}$-4, and $\rm{OOS}$-12, respectively.} 
\label{fig:example_day_dispatch}
\figspace
\end{figure}

The increase  in the number of segments $|\mathcal{S}|$ enables the electrolyzer to consume power more flexibly, as depicted in Fig.~\ref{fig:example_day_dispatch}, where the optimal power consumption schedule of the electrolyzer for one example day of the year is shown for three different numbers of segments (1, 4, and 12). It is observed that when the optimal power consumption of the electrolyzer is not constrained by wind production shortage, as on the chosen day, the optimal consumption level is always one of the piecewise linearization points. 
There are instances, e.g., hour 5 in Fig.~\ref{fig:example_day_dispatch}, where $\rm{OOS}$-1 goes into the standby state as the day-ahead price is too high for profitable hydrogen production. In contrast, $\rm{OOS}$-4 and $\rm{OOS}$-12 continue the operation in the on state, but at the power consumption level corresponding to the maximum efficiency, where hydrogen production is still profitable.

The number of segments $|\mathcal{S}|$ plays an important role in the optimal dispatch decision when the day-ahead price lies within a specific \textit{price range}. The lower bound of this price range is the price below which, for any set of segments $\mathcal{S}$, the optimal dispatch decision is always the maximum electrolyzer consumption. The upper bound is the price after which, for any set of segments $\mathcal{S}$, the electrolyzer is in standby or off state. These bounds are calculated based on the electrolyzer efficiency curve, standby power consumption, and hydrogen price, as explained in the Appendix.
Fig.~\ref{fig:price_dist} shows the distribution of day-ahead hourly prices for year 2019 in DK2 with the bounds of the price range of interest marked by the red and blue dotted lines. If the day-ahead price of a given hour lies inside this price range (green shaded in Fig.~\ref{fig:price_dist}), different dispatch decisions are taken by MILP models with different choice and number of segments.
If the day-ahead price lies outside of this range, the dispatch decision for any number of segments would be the same (i.e., produce at maximum load or cease the production) and there would be no added value of a detailed production curve.
This will be further investigated in Section \ref{sensitivity}.

\subsection{Impact of the States} 
\label{ssec:states}
We consider three cases $\rm{OOS}$, $\rm{OO}$, and $\rm{OS}$, each for both 1 and 12 segments. Recall that their corresponding MILPs are different\footnote{While we solve the proposed MILP \eqref{eq:2_objective}-\eqref{last_con3x} for $\rm{OOS}$, the MILPs presented in \cite{Online_Companion} are solved for $\rm{OO}$ and $\rm{OS}$, respectively.}.
%
%
Comparing the results of MILPs with the same number of segments, we observe $\rm{OS}$ and $\rm{OOS}$ perform almost equally,  as observed in  Fig.~\ref{fig:comparison_profit}. The reason for this is the low frequency of consecutive hours of too high day-ahead prices, where a complete shut-off would be preferred over the standby state. Over 8,760 hours, $\rm{OOS}$-1 starts up only 2 times, with a total of 286 hours offline. The difference in results obtained for $\rm{OS}$ and $\rm{OOS}$ increases if a higher standby power consumption or lower cold start-up cost for the electrolyzer is assumed, which would lead to more frequent shut-offs. On the contrary, $\rm{OO}$  earns the lowest profit, mainly due to the high start-up cost, which decreases the operational flexibility as even a short pause in production incurs a high cost. 



\begin{figure}[t]
\centering
\includegraphics[width=\linewidth]{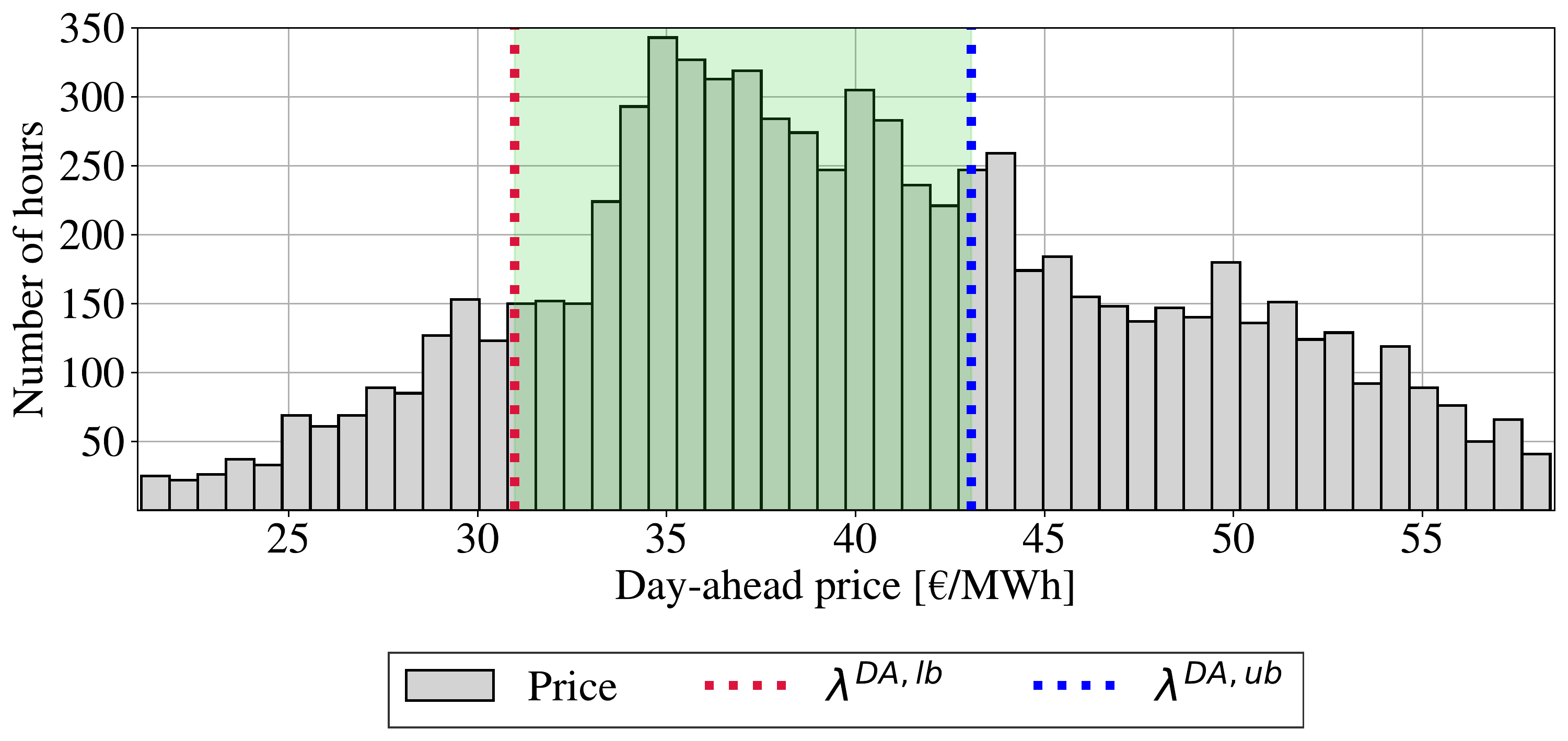}%
\caption{Histogram of day-ahead hourly electricity prices for year 2019 in DK2 (only the prices between the 5\% and 95\% quantiles are shown). Prices $\lambda^{\rm{DA, lb}}$ and $\lambda^{\rm{DA, ub}}$ (defined in the Appendix) are respectively the lower and upper bound of the price range (shaded green area) where different choices of linearization segments lead to different optimal dispatch decisions.}
\label{fig:price_dist}
\figspace
\end{figure}

\subsection{Ex-post Performance Analysis} 
\label{ssec:ex_post}
Recall that three MILPs solve the problem based on the linearized hydrogen curve.  
Through the following ex-post performance analysis, it is seen that this leads to both sub-optimal dispatch decisions and an  underestimation of the \textit{true} amount of hydrogen produced. We have already observed in   Fig.~\ref{fig:fig_eff_prod}(b) that the linearized red curve is below the original black non-linear hydrogen production curve, implying that the hydrogen production might be underestimated. This means that we can  expect to produce more hydrogen than what MILPs calculate. Such a difference is expected to be reduced by using more segments 
$|\mathcal{S}|$ to approximate the original non-linear hydrogen production curve.


Pursuing a fair comparison among models, we conduct an ex-post performance analysis. Once the MILPs are solved and the optimal power consumption $p^{\rm{e}^{*}}_t$ of the electrolyzer obtained, we re-calculate the \textit{true} amount of hydrogen produced based on the original non-linear hydrogen production curve. Note that we do not re-optimize the problem\footnote{To avoid re-optimization, we assume the extra hydrogen is directly sold to the demand and is not stored in the hydrogen storage. Otherwise, one needs to re-optimize a posteriori to optimize the operation of storage and compressor.}. We refer to the amount of extra hydrogen and its corresponding profit as 
``realized surplus". We assume that all extra hydrogen is sold at the same constant price, i.e., €2.10/kg.

\begin{figure}[t!]
\centering
\includegraphics[width=\linewidth]{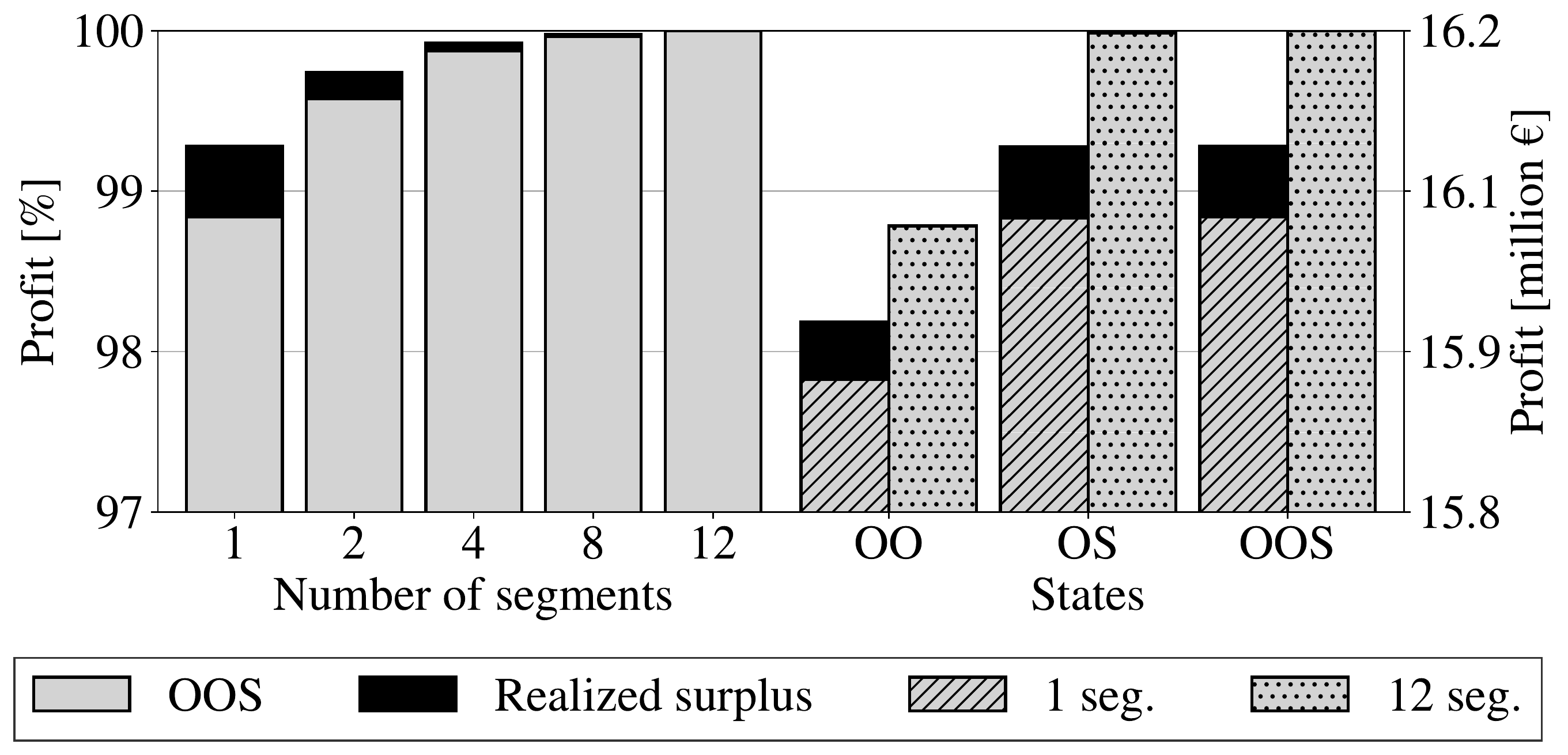}%
\caption{Estimated and realized surplus profit. The first five bars from the left correspond to $\rm{OOS}$-1 to $\rm{OOS}$-12. The next six bars show the results for $\rm{OO}$-1, $\rm{OO}$-12, $\rm{OS}$-1, $\rm{OS}$-12, $\rm{OOS}$-1, and $\rm{OOS}$-12, respectively. The right vertical axis is the profit in million €, whereas the left vertical axis is the relative profit in \% in comparison to the highest profit achieved by $\rm{OOS}$-12.} 
\label{fig:comparison_profit}
\figspace
\end{figure}

\begin{figure}[t!]
\centering
\includegraphics[width=\linewidth]{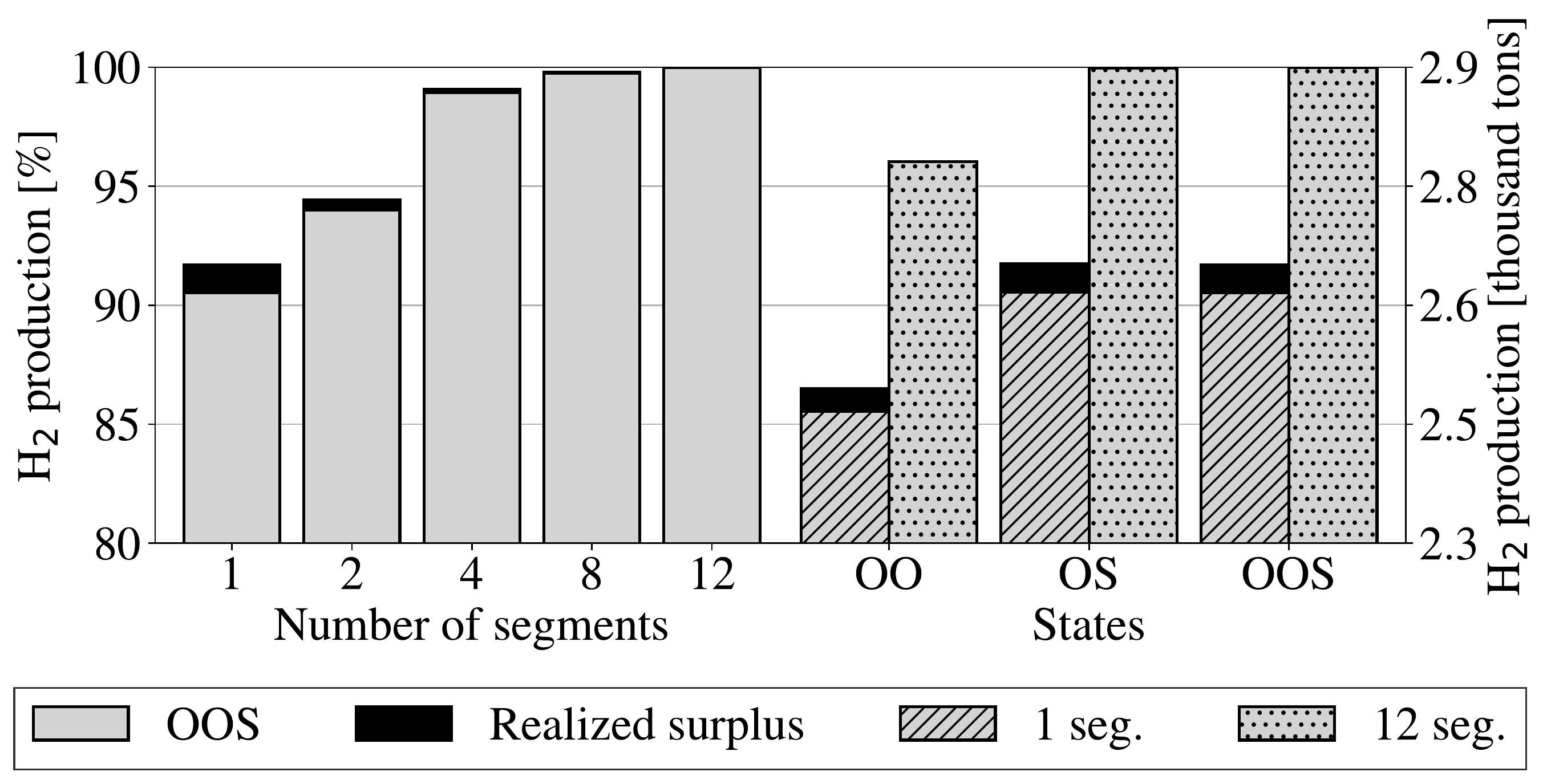}%
\caption{Estimated and realized surplus hydrogen produced. The first five bars from the left correspond to $\rm{OOS}$-1 to $\rm{OOS}$-12. The next six bars show the results for $\rm{OO}$-1, $\rm{OO}$-12, $\rm{OS}$-1, $\rm{OS}$-12, $\rm{OOS}$-1, and $\rm{OOS}$-12, respectively.}
\label{fig:comparison_prod}
\figspace
\end{figure}

Fig.~\ref{fig:comparison_profit} provides the \textit{estimated} and \textit{realized surplus} profit among different cases. The estimated profit (gray area) is the optimal value obtained for the objective function of the corresponding MILP, while the realized profit (dark area), calculated ex-post, takes into account the profit of selling extra hydrogen. Similarly,  Fig.~\ref{fig:comparison_prod} shows the total estimated and realized surplus hydrogen produced. Note that the compressor would need to consume more power ($\sim$ 1 MWh/ton) due to extra hydrogen.
We draw two conclusions from Fig.~\ref{fig:comparison_profit} and \ref{fig:comparison_prod}:




(1) \textit{Realized surplus}: This surplus for profit and hydrogen production is reduced by increasing the number of segments, due to the improved approximation of the original non-linear curve. The realized surplus profit decreases from €71,199 (0.44\%) for $\rm{OOS}$-1  to €602 (below 0.01\%) for $\rm{OOS}$-12. Similarly, the hydrogen production surplus is significantly decreased, yielding a realized surplus of $\sim 34$ tons (1.27\%) for $\rm{OOS}$-1 and only 0.3 tons (0.01\%) for $\rm{OOS}$-12. By choosing a low number of segments, the hydrogen production is underestimated which may lead to logistic issues and inefficiencies in the real-life operation of the hybrid power plant.

(2) \textit{Ex-post profit and hydrogen production}: Adding more electrolyzer details (segments or/and states) always leads to an increase in the ex-post profit. To compare various models, $\rm{OOS}$-12 is taken as a benchmark, as it leads to the highest profit. First, the impact of the number of segments is examined, while keeping the number of states fixed and equal to 3. The ex-post profit reduction applying 1 instead of 12 segments is 0.72\%, corresponding to around 117.6 k€ for the entire hybrid power plant. The ex-post hydrogen production is increased by 8.32\%, corresponding to around  241 tons. This percentage deviation is notably higher in part because the increase in hydrogen profit is dampened by the reduction in electricity profit (3.86\% electricity profit increase for 1 segment compared to 12 segments). For $\rm{OOS}$-1, the profit share of selling hydrogen is much lower than the profit share of selling electricity (around 34\%). By introducing more segments, the contribution of hydrogen sales is increased to 38\% at the expense of electricity sales. More profit and different business models are therefore unlocked by including more electrolyzer details in the MILP formulation. Fig.~ \ref{fig:comparison_profit} and \ref{fig:comparison_prod} show that the errors are considerably reduced by implementing 4 segments instead of 1.
Second, we assess the impact of the states on the ex-post profit and hydrogen production. While  $\rm{OS}$ performs just as well as  $\rm{OOS}$ as described in Section \ref{ssec:states}, $\rm{OO}$ with 12 segments results in a 1.22\% lower ex-post profit, and in a 4\% lower hydrogen production. For $\rm{OO}$-1, a profit reduction of around 1.8\% and a reduced hydrogen production of 13.5\% are observed, compared to the benchmark. Finally, we observe that neglecting the standby state in the model formulation leads to the worst outcome in terms of profit and hydrogen production potential.

\begin{table}[t]
    \caption{Computational aspects}
    \label{tab:comp_aspects}
    \centering
    \begin{tabular}{lrr}
    \toprule
        Case  & Computational time {[}s{]} & No. of binary variables \\

    \midrule
    $\rm{OS}$-1  & 1.4 & 2×8760 \\
    $\rm{OS}$-12  & 12.7 & 13×8760 \\ \hline
    $\rm{OOS}$-1  & 137.8 & 5×8760 \\
    $\rm{OOS}$-2  & 135.8 & 6×8760 \\
    $\rm{OOS}$-4 & 236.3 & 8×8760 \\
    $\rm{OOS}$-8  & 350.3 & 12×8760 \\
    $\rm{OOS}$-12  & 473.7 & 16×8760 \\ \hline
    $\rm{OO}$-1  & 767.1 & 3×8760 \\
    $\rm{OO}$-12  & 1,763.1 & 14×8760 \\
\bottomrule
\end{tabular}
\figspace
\end{table}
\subsection{Computational Analysis} 
\label{ssec:comp}
All MILPs have been solved using the Gurobi solver in Julia on a MacBook Pro M1 2020 with 16 GB RAM.
The optimality gap is fixed to $0.01$\% when we solve every MILP. 
The increase in the number of linearization segments $|\mathcal{S}|$ leads to an increase in computational time due to introducing more binary variables. For $\rm{OOS}$, the computational time is increased from  138 seconds for 1 segment to 474 seconds for 12 segments, as reported in  Fig.~\ref{tab:comp_aspects}. Removing the off state significantly reduces the computational time, with $\rm{OS}$-1 being by far the fastest MILP to be solved (1.4 seconds). The $\rm{OO}$ models require the highest computational time, although they embody fewer binary variables than their $\rm{OS}$ and $\rm{OOS}$ counterparts. We hypothesize the reason is that the start-up cost constraints with inter-temporal nature are more often active when the option of standby state is not present. Therefore, we do not recommend using $\rm{OO}$ as its corresponding profit is the lowest among all cases (Fig.~\ref{fig:comparison_profit}), and it is being solved comparatively slower. Further, if computational efficiency is crucial, it may be beneficial to neglect the off state and run the $\rm{OS}$ model for improved computational performances. In general, the computational time increases with the number of segments but is deemed reasonable for the $\rm{OS}$ and $\rm{OOS}$ models, considering that our optimization problem is run over 8,760 hours. As operational problems are typically solved for a shorter time horizon, e.g., 24 hours for day-ahead scheduling, the computational cost of adding more details to the electrolyzer would be minimal.


\subsection{Sensitivity Analysis with Respect to Input Data}
\label{sensitivity}
In the previous sections, we have shown that adopting a simplified electrolyzer model can lead to an underestimation of the profit and hydrogen production for the hybrid power plant. We have also shown that the benefit of added details is case-specific, and depends on the input parameters. We now aim at assessing the impact of input parameters and system configuration on these results, through a sensitivity analysis. In particular, we will focus on wind over electrolyzer capacity ratio, hydrogen demand over electrolyzer capacity ratio, and the hydrogen price. The sensitivity analysis is performed on the $\rm{OOS}$-1  and $\rm{OOS}$-12  models.

\subsubsection{Wind size}
Recall from  Fig.~\ref{tab:case_details} that the wind farm capacity is 2 times that of the electrolyzer. To assess the impact of the wind-to-electrolyzer capacity ratio, two additional cases are considered, under which such a ratio is 1, 2 (reference), and 8. When this ratio is reduced from 2 to 1, the number of hours where the power input to the electrolyzer is limited by the wind availability is increased from 5,326 to all hours. Conversely, when the ratio is increased from 2 to 8, the number of power-limited hours is reduced to 1,236.
We observe that the realized surplus for hydrogen production increases with the number of hours with limited wind power. The reason for this is that the piecewise approximation is exact only on the linearization points, and the limited wind availability forces the electrolyzer to operate out of those points. Conversely, when the number of wind power-limited hours is reduced, the electrolyzer operates more often on the linearization points, where the approximation is exact. It follows that the underestimation of hydrogen production is greater the more the electrolyzer is limited from operating at the linearization points.  With a wind-to-electrolyzer ratio of 1, the difference in ex-post hydrogen production between 1 and 12 segments is 13\%, which is reduced to 3\% when the ratio increases to 8. Therefore, incorporating electrolyzer details is crucial for hybrid power plants where the wind-to-electrolyzer capacity ratio is small.

\subsubsection{Hydrogen demand size}
To investigate the sensitivity of optimization outcomes with respect to the hydrogen demand, the minimum daily demand is doubled, corresponding to around 8 full-load hours of hydrogen production.  We observe that the impact of adding more segments to the electrolyzer production curve diminishes when the demand constraint is tighter, i.e., with a higher minimum daily demand. For the case with the reference demand, the difference between the ex-post profit for $\rm{OOS}$-12 and $\rm{OOS}$-1  is 8\%.  This difference, when the hydrogen demand is doubled, is reduced to 2\%. The increase in demand forces the electrolyzer to operate more frequently at its maximum load, where both $\rm{OOS}$-1 and $\rm{OOS}$-12 share the same linearization point and efficiency.
\subsubsection{Hydrogen price}
To explore the impact of the hydrogen price, we increase it from €2.1/kg to €5/kg. 
As already discussed in Section \ref{ssec:segments}, adding more segments impact the optimal solution and profit as long as the electricity price in the given hour is in the range $[\lambda^{\rm{DA, lb}}, \lambda^{\rm{DA, ub}}]$, shown in  Fig.~\ref{fig:price_dist}. Since $\lambda^{\rm{DA, lb}}$ and $\lambda^{\rm{DA, ub}}$ are proportional to the hydrogen price, by increasing the hydrogen price, the range is moved towards higher electricity prices, where the frequency of occurrence is reduced. When the MILP is solved with the hydrogen price of €5/kg, it is more frequently optimal to operate the electrolyzer at full load (39\% of the time, compared to 11\% for the case with the hydrogen price of €2.1/kg) and the linearization segments are utilized less. This also results in a significantly decreased computational time (below 20 seconds for $\rm{OOS}$-12). The profit contribution from the hydrogen sale is increased significantly to 92\%. The ex-post profit and hydrogen production difference between $\rm{OOS}$-1 and $\rm{OOS}$-12 are reduced to 0.01\% and 0.03\%, respectively (they are 0.72\% and 8.32\% for the €2.1/kg case).

The modeling of segments is relevant if higher hydrogen prices are coupled with also higher electricity prices. In this way, the electricity price range $[\lambda^{\rm{DA, lb}}, \lambda^{\rm{DA, ub}}]$ would still be overlapping with the majority of day-ahead price occurrences. For example, we test an artificial  case where the day-ahead electricity price time series was multiplied by a constant factor to increase the mean price to around €90/MWh (similar to the mean value for 2021 in DK2). In this case, with the hydrogen price of €5/kg, similar results to the 2019 test case with the hydrogen price of €2.1/kg were obtained in terms of the impact of the number of segments. For a given hydrogen price and efficiency curve, checking if the price range $[\lambda^{\rm{DA, lb}}, \lambda^{\rm{DA, ub}}]$ overlaps with the expected electricity price is therefore crucial to assess a priori the impact of choosing a simplified model for the production curve (e.g., 1 linearization segment only) and support the modeling choices.

\section{Discussion and Conclusion}

\label{conclusion}
Several studies have focused on the optimal dispatch of hybrid renewable-hydrogen power plants assuming simplified models for the electrolyzer component. This paper investigates the impact of choosing different levels of operational details for the electrolyzer model on the dispatch decisions, profit, the amount of hydrogen produced, and computational time. The impact of two modeling choices is considered: the operating states (on, off, standby), and the number of segments used to linearize the hydrogen production curve. The problems are formulated as MILPs, where the number of binary variables depends on the number of states and segments.

For fixed states, adding more linearization segments for approximating the hydrogen production curve results in a higher profit, and a reduced surplus in the ex-post profit calculation, meaning that the model is able to estimate the actual cost and revenue streams more accurately. Moreover, a better estimation of the produced hydrogen is achieved. In fact, the linearization results in an underestimation of the produced hydrogen, but the underestimation is reduced by increasing the number of segments. Apart from introducing errors in the actual realized profit, thus potentially impacting the investment decisions in these types of technologies, the systematic underestimation of the hydrogen produced by the electrolyzer might introduce logistical  inefficiencies, e.g., truck scheduling, and storage discharging/filling. 

The impact of adding more piecewise segments to the hydrogen production curve depends on the distribution of day-ahead electricity prices in the given time horizon. The model formulations with 1 and 12 segments take significantly different dispatch decisions when the day-ahead electricity price is within a certain range, which depends on the electrolyzer efficiency (minimum and maximum) and the hydrogen price. Out of this day-ahead electricity price range, the model with 1 and 12 segments takes the same dispatch decisions. Therefore, the value of adding more details to the hydrogen production curve could differ by varying input data and case studies. It is observed that this value decreases when the electrolyzer operates less at partial loading, e.g. when the input power is less limited by available wind power or with high-demand constraints. In this paper, revenues from other than the day-ahead market are not considered but this may also impact the dispatch strategy and therefore benefit from more segments.

Choosing to represent only on and off states leads to the highest profit underestimation and worst ex-post performance while modeling only on and standby states lead to similar profit and dispatch decisions to the three-state model. This result is, however, significantly affected by the assumption made on the standby power consumption of the electrolyzer and its start-up cost. These parameters are highly uncertain due to the lack of data on large-scale electrolyzers. 

In conclusion, adopting more simplified models for the electrolyzer always leads to a reduced profit and sub-optimal scheduling. However, the impact of adding more details may vary depending on the case study considered and especially the range of day-ahead electricity prices, hydrogen price, wind power production compared to the electrolyzer installed capacity, standby power consumption, and  start-up cost. Among all considered models, the most complete one (three states with 12 segments) was solved for a 1-year horizon in less than 10 minutes. The increase in computational time by adding more details would be marginal if a day-ahead scheduling problem is considered instead. Moreover, reducing the three-state model to two states only is not always faster, as it was observed that the two-state on-off model with 12 segments was the longest to solve among all the cases considered.  A more detailed representation of the electrolyzers should be preferred for operational problems. For investment problems, we hypothesize that it may be adequate to adopt a more simplified model of the electrolyzer, but this should be further assessed and it was out of the scope of the current paper.

Further research should be conducted to assess the impact of modeling choices when additional revenue streams are considered, such as flexibility provisions in ancillary service markets, which may impact the dispatch decisions of the hybrid power plant. 


\begin{appendix}
\label{AppA}
\setcounter{equation}{0}
\renewcommand{\theequation}{A\thesubsection.\arabic{equation}}

This appendix provides an analytical formulation for the day-ahead electricity price range bounds ($\lambda^{\rm{DA,lb}}$ and $\lambda^{\rm{DA,ub}}$ in Fig.~\ref{fig:price_dist}), where different dispatch decisions are made by MILP models based on the choice of linearization segments. 
To derive these bounds we assume that the electrolyzer operation is not constrained by limited wind availability or demand constraints.
We also neglect the compressor power consumption. Therefore, for a given hour, the power sold to the grid is calculated as the wind power production $P^{\rm{w}}$ minus the electrolyzer power consumption.

\begin{figure}[t]
\centering
\includegraphics[width=0.76\linewidth]{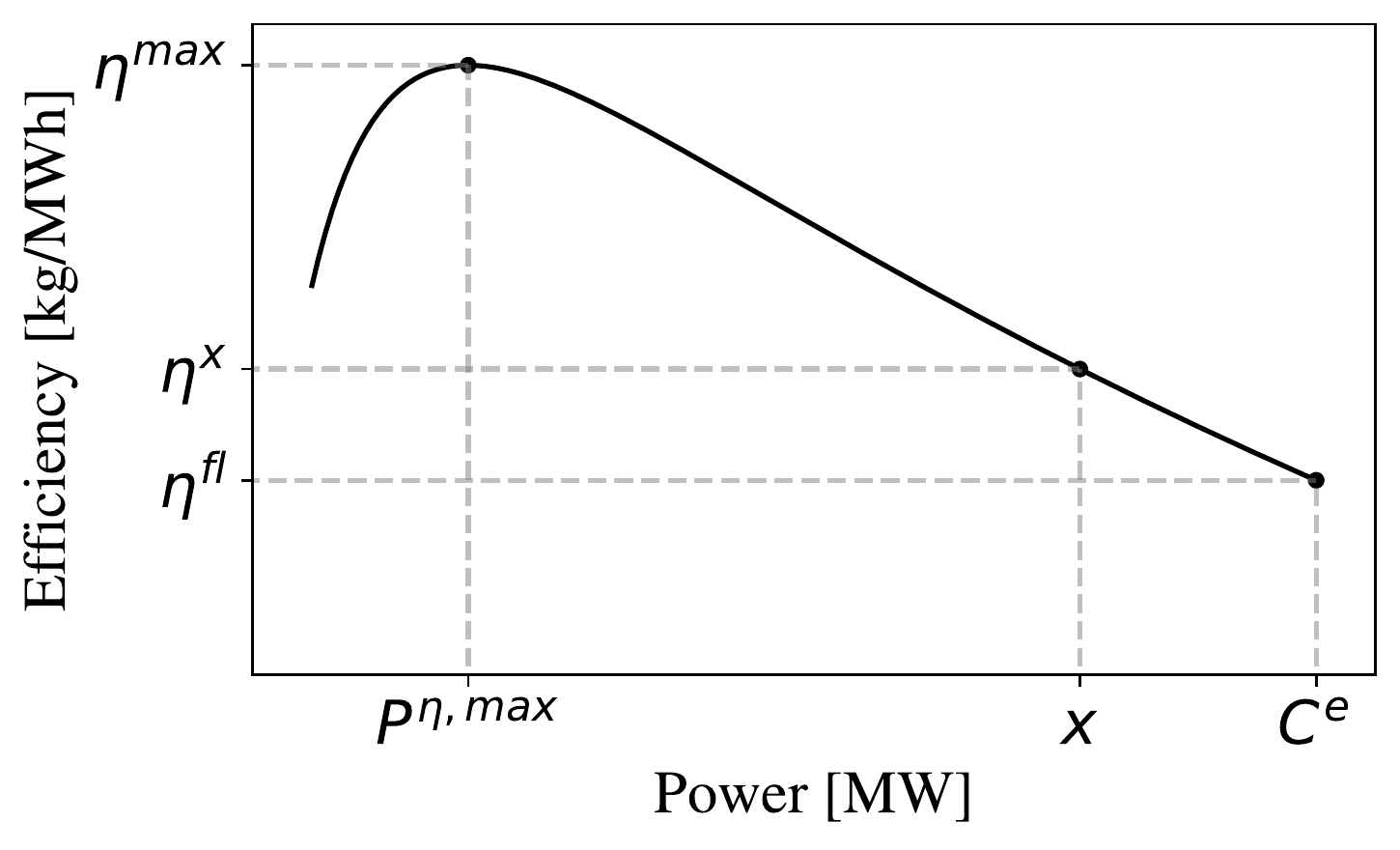}%
\vspace{-0.15cm}
\caption{Efficiency curve with relevant points to calculate the bounds of the electricity price range wherein electrolyzer operational details matter.}
\label{fig:efficiency_appendix}
\figspace
\end{figure}

The upper bound $\lambda^{\rm{DA,ub}}$ 
refers to a price for which it becomes more profitable for the hybrid power plant to operate the electrolyzer at point $P^{\rm{\eta,max}}$ with efficiency $\eta^{\rm{max}}$ (see Fig.~\ref{fig:efficiency_appendix}) than to keep it in the standby mode, i.e.,

{\small
\begin{align}
    &\lambda^{\rm{DA}}(P^{\rm{w}}-P^{\rm{\eta,max}})+\lambda^{\rm{h}}\eta^{\rm{max}}P^{\rm{\eta,max}} > \lambda^{\rm{DA}}(P^{\rm{w}}-P^{\rm{sb}}),\label{eq:upper_bound1}
\end{align}
}%
where the left-hand side is the profit of the hybrid power plant when the electrolyzer consumes $P^{\rm{\eta,max}}$, and the right-hand side is that when the electrolyzer is in the stand-by mode. This gives the upper bound $\lambda^{\rm{DA,ub}}$ as

{\small
\vspace{-0.15cm}
\begin{align}
    &\lambda^{\rm{DA,ub}} =  \frac{\lambda^{\rm{h}} \eta^{\rm{max}}P^{\rm{\eta,max}} }{P^{\rm{\eta,max}}-P^{\rm{sb}}}.\label{eq:upper_bound}
\end{align}
}%
For prices higher than $\lambda^{\rm{DA,ub}}$, for any segment choice, it is not profitable for the hybrid power plant to produce hydrogen. 

The lower bound $\lambda^{\rm{DA,lb}}$ is the price below which, for any choice of segments, it is always optimal to operate the electrolyzer at full load $C^{\rm{e}}$ with efficiency $\eta^{\rm{fl}}$. Consider the closest linearization point to the full load power, which is hypothetically denoted by $x$ in Fig.~\ref{fig:efficiency_appendix}. Depending on $\lambda^{\rm{DA}}$, it could be profitable to reduce the electrolyzer power consumption to operate at a higher efficiency $\eta^{\rm{x}}>\eta^{\rm{fl}}$, i.e., 

{\small
\vspace{-0.15cm}
\begin{align}
    & \lambda^{\rm{DA}}(P^w-C^{\rm{e}})+\lambda^{\rm{h}}\eta^{\rm{fl}}C^{\rm{e}} > \lambda^{\rm{DA}}(P^w-x)+\lambda^{\rm{h}}\eta^{\rm{x}}x. \label{eq:lower_bound_equation}
\end{align}
}%
The lowest $\lambda^{\rm{DA}}$ for which \eqref{eq:lower_bound_equation} fulfils is obtained as $x$ tends to $C^{\rm{e}}$ and it can therefore be expressed as the following limit:

{\small
\vspace{-0.3cm}
\begin{align}
        &\lambda^{\rm{DA,lb}} = \lim_{x\to C^{\rm{e}}} \frac{\lambda^{\rm{h}} \eta^{\rm{fl}}C^{\rm{e}}-\lambda^{\rm{h}} \eta^{\rm{x}}x}{C^{\rm{e}}-x} ,\label{eq:lower_bound}
\end{align}
}%

The limit in \eqref{eq:lower_bound} is an indeterminate form that can be solved using L'Hôpital's rule of limit theory, obtaining

{\small
\vspace{-0.3cm}
\begin{align}
     \lambda^{\rm{DA,lb}}   = \lambda^{\rm{h}}(\eta^{\rm{fl}}+C^{\rm{e}} \eta'(x)|_{x=C^{\rm{e}}}), \label{eq:lower_bound_limit}
\end{align}
}%
where $\eta'(x)|_{x=C^{\rm{e}}}$ is the derivative of the efficiency curve calculated at full load power. Since there is no analytical formulation for the efficiency curve, $\eta'(x)|_{x=C^{\rm{e}}}$ can be estimated using finite differences.

\end{appendix}

\section*{Acknowledgement }
This research was supported by the Energy Cluster Denmark through the ``Sustainable P2X Business Model" project, and by the Danish Energy Development Programme (EUDP) through the HOMEY project (64021-7010). We would like to thank  Jens Jakob Sørensen (Ørsted), Alexander Holm Kiilerich (Ørsted), Roar Hestbek Nicolaisen (Hybrid Greentech), Yannick Werner (DTU), and Matěj Novotný  for collaborations, thoughtful discussions, and constructive feedback. 
\bibliography{8_References.bib} 
\bibliographystyle{ieeetr}



\end{document}